\documentclass[10pt]{amsart}

\include{code_preamble.tex}

\usepackage{amsmath}
\usepackage{amssymb}
\usepackage{amsthm}
\usepackage{hyperref}
\usepackage{graphicx}
\usepackage{cleveref}
\usepackage{xcolor}
\usepackage{wrapfig}
\usepackage[style=alphabetic]{biblatex}

\numberwithin{equation}{section}

\theoremstyle{definition}

\newtheorem*{remark*}{Remark}

\DeclareMathOperator{\im}{im}

\newcommand{\C}{\mathbb{C}}

\title{SeqSee: A schema-based approach to spectral sequence visualization}
\author{Joey Beauvais-Feisthauer}
\email{joeybf@wayne.edu}
\author{Dan Isaksen}
\email{isaksen@wayne.edu}
\address{Department of Mathematics, Wayne State University, Detroit, Michigan 48009, USA}

\date{30 January 2025}

\thanks{The authors were partially supported by National Science Foundation Grant DMS-2202267.}



\begin{document}

\begin{abstract}
	We present \texttt{SeqSee}, a software system that addresses spectral sequence visualization through a schema-based approach. By introducing a standardized JSON schema as an intermediate representation, \texttt{SeqSee} decouples the mathematical computations of spectral sequences from their visualizations.  We demonstrate the system through a case study of the classical and $\C$-motivic Adams spectral sequences.
\end{abstract}

\maketitle

\section{Introduction}

Computational stable homotopy theory is currently undergoing a revolution using machines to generate, organize, and interpret data.  The fundamental problem is to compute the stable homotopy ring, which is a very deeply structured object \cite{MR4250190}.

Spectral sequences are fundamental tools that help us compute through successive approximations \cite{MR1793722}.  Effective use of spectral sequences typically requires visualizations to make the data accessible.  We aim to make it easy for mathematicians to produce visualizations of their spectral sequences, so they can focus their efforts on generating the underlying mathematical data.

\texttt{SeqSee} (pronounced ``seek-see'') is a generic visualization
tool for spectral sequence data.  We present a JSON schema that serves as
an intermediate representation, midway between raw mathematical data and explicit code that produces graphics, that serves as a clear boundary between mathematical computation and visual presentation.

\texttt{SeqSee} takes a JSON file as input, conforming to the \texttt{SeqSee}
schema, and outputs a self-contained HTML file suitable for display in any browser.
This file includes an
SVG figure representing the spectral sequence along with JavaScript for
interactivity.

\texttt{SeqSee} is general enough to handle all of the many spectral sequences that the second author has studied in the last 20 years.  These include spectral sequences of a purely algebraic nature, as well as those with homotopical content.   They include spectral sequences that are associated with
classical stable homotopy,
$\C$-motivic stable homotopy \cite{isaksen_2022_6547197, MR2591921, isaksen_2022_6987157, isaksen_2022_6987227, MR4046815, MR4806410},
$\mathbb{R}$-motivic homotopy \cite{MR4461846, MR3572357, MR4777709},
and $C_2$-equivariant stable homotopy \cite{MR4041284, guillou2024c2equivariantstablestems}.
\texttt{SeqSee} unifies the display of spectral sequences from all of these projects, as well as similar types of projects of the future.

\begin{figure}
	\includegraphics[width=0.7\linewidth,trim={1.6cm 14.7cm 8.3cm 1cm},clip]{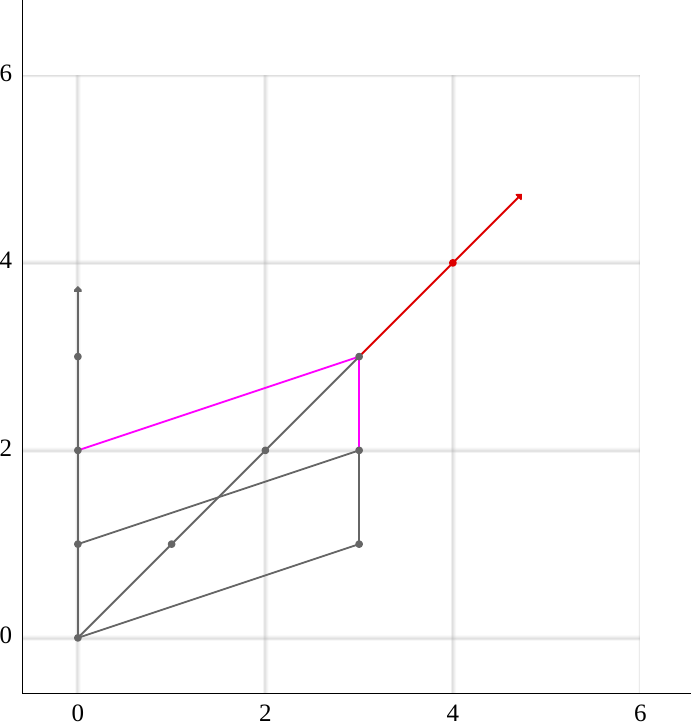}
	\caption{\texttt{SeqSee} visualization of the initial part of the $\C$-motivic Adams spectral sequence, displaying coordinate axes, nodes, and edges.}
\end{figure}

Currently, there are at least three entirely separate machined-based projects
\cite{CCBF, Bruner, Lin}
for producing large-scale data that is relevant to the Adams spectral sequence.
The authors foresee a need for these, and other more loosely related projects, to interoperate.  \texttt{SeqSee} aims to unite at least the visualization of spectral sequences
that have been computed with different platforms, even though their underlying data are
stored in distinct formats.

Other software
packages for spectral sequence computation \cite{CCBF, Lin} provide
similar functionality for spectral sequence visualization.  \texttt{SeqSee} aims
to be more expressive and more user-friendly than these other options. Also, while those other packages provide more features, they make use of linear algebra for introspection. \texttt{SeqSee} operates purely at the graphical level, which we believe to be an advantage.

\section{Mathematical Foundation}
A spectral sequence is a sequence $(E_r,d_r)_{r\geq 1}$ of differential graded modules  where each page $E_{r+1}$ is the homology of the previous page $E_r$ with respect to the differential $d_{r}$. The equation
\[
	E_{r+1} = H(E_r,d_r) = \frac{\ker d_r}{\im d_r}.
\]
describes this structure.
Spectral sequences are inherently bigraded.
The differential $d_r$ operates between the homogeneous components
of $E_r$ with a specific bidegree shift that depends on $r$ and also on the structure
and choice of conventions for the spectral sequence at hand.
\texttt{SeqSee} is agnostic to the choice of grading conventions on the differentials.

Every spectral sequence can be summarized in a bigraded visual chart.
A spectral sequence chart displays information about different bigradings in different locations across a 2-dimensional plane.

Such charts do not attempt to capture the full structure of the dataset, which typically is vastly more rich than what could possibly be contained in a visual chart. Nevertheless, charts are essential
for studying the overall structure of a spectral sequence,
and they are likely to continue in this role.  

From the perspective of visualization, the essential properties of a spectral sequence chart is that it is a bigraded coordinate system, indexed on integers.  Moreover, there are certain elements that called \emph{nodes} at certain coordinates, and there are certain relations amongst pairs of nodes called \emph{edges}.  In simplest terms, \texttt{SeqSee} plots dots and lines at locations that the user specifies.

\section{Case study: The classical and $\C$-motivic Adams spectral sequences}

The capabilities of \texttt{SeqSee} are particularly well demonstrated through its handling of the classical and $\C$-motivic Adams spectral sequences \cite{isaksen_2022_6987157}. These spectral sequences, which compute stable homotopy groups of spheres, provide excellent test cases due to their complexity and importance in algebraic topology.  The \textsc{Figure} shows a very small
sample of the $\C$-motivic Adams spectral sequence.
In this example, nodes represent a basis for the spectral sequence, and edges represent both differentials and also some of the multiplicative structure.

The visualization of the Adams spectral sequence highlights several key features of \texttt{SeqSee}. The system efficiently handles the dense network of relationships between nodes, maintains precise positioning of nodes in the bigraded structure, and provides clear visual representation of differentials and extensions.

The \texttt{SeqSee} GitHub repository \cite{Beauvais-Feisthauer_SeqSee_2025} includes the raw data
for the case study in CSV format.
The repository also includes a tool called \texttt{jsonmaker},
which can convert the CSV files
into JSON files following our schema.
This script is specific to the format of those CSV
files, but serves as a template for creating tools to generate JSON data
that \texttt{SeqSee} can display. It is expected that any software that
outputs spectral sequence data will require a customized tool to produce
JSON files that follow the \texttt{SeqSee} schema.
Finally, the repository includes the HTML output.

\section{Schema Design}

The core innovation of \texttt{SeqSee} lies in its JSON schema, which serves as a specialized intermediate representation for spectral sequences. This schema preserves just enough mathematical structure for the needs of visualization. Here is a simplified example of the schema structure:

\begin{Verbatim}[commandchars=\\\{\}]
  \PY{p}{\PYZob{}}
  \PY{+w}{  }\PY{n+nt}{\PYZdq{}header\PYZdq{}}\PY{p}{:}\PY{+w}{ }\PY{p}{\PYZob{}}
  \PY{+w}{    }\PY{n+nt}{\PYZdq{}metadata\PYZdq{}}\PY{p}{:}\PY{+w}{ }\PY{p}{\PYZob{}}\PY{n+nt}{\PYZdq{}title\PYZdq{}}\PY{p}{:}\PY{+w}{ }\PY{l+s+s2}{\PYZdq{}Example Spectral Sequence\PYZdq{}}\PY{p}{\PYZcb{},}
  \PY{+w}{    }\PY{n+nt}{\PYZdq{}chart\PYZdq{}}\PY{p}{:}\PY{+w}{ }\PY{p}{\PYZob{}}
  \PY{+w}{      }\PY{n+nt}{\PYZdq{}width\PYZdq{}}\PY{p}{:}\PY{+w}{ }\PY{p}{\PYZob{}}\PY{n+nt}{\PYZdq{}min\PYZdq{}}\PY{p}{:}\PY{+w}{ }\PY{l+m+mi}{0}\PY{p}{,}\PY{+w}{ }\PY{n+nt}{\PYZdq{}max\PYZdq{}}\PY{p}{:}\PY{+w}{ }\PY{l+m+mi}{5}\PY{p}{\PYZcb{},}
  \PY{+w}{      }\PY{n+nt}{\PYZdq{}height\PYZdq{}}\PY{p}{:}\PY{+w}{ }\PY{p}{\PYZob{}}\PY{n+nt}{\PYZdq{}min\PYZdq{}}\PY{p}{:}\PY{+w}{ }\PY{l+m+mi}{0}\PY{p}{,}\PY{+w}{ }\PY{n+nt}{\PYZdq{}max\PYZdq{}}\PY{p}{:}\PY{+w}{ }\PY{l+m+mi}{5}\PY{p}{\PYZcb{}}
  \PY{+w}{    }\PY{p}{\PYZcb{}}
  \PY{+w}{  }\PY{p}{\PYZcb{},}
  \PY{+w}{  }\PY{n+nt}{\PYZdq{}nodes\PYZdq{}}\PY{p}{:}\PY{+w}{ }\PY{p}{\PYZob{}}
  \PY{+w}{    }\PY{n+nt}{\PYZdq{}1\PYZdq{}}\PY{+w}{ }\PY{p}{:}\PY{+w}{ }\PY{p}{\PYZob{}}\PY{n+nt}{\PYZdq{}x\PYZdq{}}\PY{p}{:}\PY{+w}{ }\PY{l+m+mi}{0}\PY{p}{,}\PY{+w}{ }\PY{n+nt}{\PYZdq{}y\PYZdq{}}\PY{p}{:}\PY{+w}{ }\PY{l+m+mi}{0}\PY{p}{,}\PY{+w}{ }\PY{n+nt}{\PYZdq{}label\PYZdq{}}\PY{p}{:}\PY{+w}{ }\PY{l+s+s2}{\PYZdq{}\PYZdl{}1\PYZdl{}\PYZdq{}}\PY{+w}{  }\PY{p}{\PYZcb{},}
  \PY{+w}{    }\PY{n+nt}{\PYZdq{}h0\PYZdq{}}\PY{p}{:}\PY{+w}{ }\PY{p}{\PYZob{}}\PY{n+nt}{\PYZdq{}x\PYZdq{}}\PY{p}{:}\PY{+w}{ }\PY{l+m+mi}{0}\PY{p}{,}\PY{+w}{ }\PY{n+nt}{\PYZdq{}y\PYZdq{}}\PY{p}{:}\PY{+w}{ }\PY{l+m+mi}{1}\PY{p}{,}\PY{+w}{ }\PY{n+nt}{\PYZdq{}label\PYZdq{}}\PY{p}{:}\PY{+w}{ }\PY{l+s+s2}{\PYZdq{}\PYZdl{}h\PYZus{}0\PYZdl{}\PYZdq{}}\PY{p}{\PYZcb{}}
  \PY{+w}{  }\PY{p}{\PYZcb{},}
  \PY{+w}{  }\PY{n+nt}{\PYZdq{}edges\PYZdq{}}\PY{p}{:}\PY{+w}{ }\PY{p}{[}\PY{+w}{ }\PY{p}{\PYZob{}}\PY{n+nt}{\PYZdq{}source\PYZdq{}}\PY{p}{:}\PY{+w}{ }\PY{l+s+s2}{\PYZdq{}1\PYZdq{}}\PY{p}{,}\PY{+w}{ }\PY{n+nt}{\PYZdq{}target\PYZdq{}}\PY{p}{:}\PY{+w}{ }\PY{l+s+s2}{\PYZdq{}h0\PYZdq{}}\PY{p}{\PYZcb{}}\PY{+w}{ }\PY{p}{]}
  \PY{p}{\PYZcb{}}
  \end{Verbatim}

The schema is organized into three main components, each serving a distinct purpose in the representation of spectral sequences. The header section contains metadata and global configuration settings. The nodes section represents individual elements of the spectral sequence, encoding their position in the bigraded structure along with mathematical properties and visual attributes. The edges section captures the relationships between nodes, such as differentials and extensions, encoding their mathematical properties and visual attributes.

\section{Implementation}

The implementation of \texttt{SeqSee} follows a modular architecture that reflects the separation between mathematical structure and visual representation. A Python-based processing pipeline transforms spectral sequence data through two stages. The process begins with the conversion of raw data, often in a bespoke format stored in CSV files or in databases, into the standardized JSON representation.

\texttt{SeqSee} users bear responsibility for converting their data into JSON.  Different projects store their data in fundamentally different ways, so there can be no all-purpose tool.  However, we provide a sample tool \texttt{jsonmaker}, customized for one particular project \cite{isaksen_2022_6987157}, that other users can adapt for their own purposes.

Once data is in the standardized format, \texttt{SeqSee} generates interactive HTML. The visualization layer uses SVG for vector graphics, ensuring crisp rendering at any scale, while KaTeX handles the display of text and mathematical notation.


The source code of \texttt{SeqSee} is MIT-licensed and is available on GitHub at \cite{Beauvais-Feisthauer_SeqSee_2025}. The repository includes the full JSON schema specification, the Python implementation of the processing pipeline, and a collection of example datasets, including the classical and $\C$-motivic Adams spectral sequences discussed above in the Case Study. The codebase follows best practices for scientific software development, including thorough documentation.

A guiding principle is to keep the technical barrier to entry as low as possible.  Familiarity with Python is required.
\texttt{SeqSee} is packaged with \texttt{uv} for easy packaging and development.

We have attempted to minimize dependencies, but
the HTML files generated by \texttt{SeqSee} do have a few dependencies and
cannot be used offline.

\section{Future outlook}

\texttt{SeqSee} is designed with flexibility in mind, making it
straightforward to extend its capabilities.
The project genuinely welcomes contributions from the community.
We encourage collaboration on extensions to \texttt{SeqSee} to satisfy specific needs.
Some short-term possible improvements are:
\begin{itemize}
	\item
	      \textbf{Customizable Styles}: Add support for more node and edge styles, shadow effects, or different node shapes.
	\item
	      \textbf{Advanced Interactivity}: Implement features like node/edge selection, more hover effects, and dynamic highlighting.
	\item
	      \textbf{Improved Coordinate Systems}: Allow more customization of axes, gridlines, and coordinates for specialized use cases.
	\item
	      \textbf{Enhanced Output Formats}: Support additional export formats,
	      such as PDFs, Tikz figures, or raw SVG files for presentations and  publications.
\end{itemize}

We foresee further uses of \texttt{SeqSee} beyond stable homotopy theory.  For example, \texttt{SeqSee} is well-suited for
equivariant slice spectral sequences \cite{Chatham} and
unstable homotopy group computations \cite{MR1232201}.

\begin{remark*}
Since the first version of this article was posted, there have been several updates to \texttt{SeqSee}.  Highlights include:
\begin{itemize}
\item
customizable display styles, such as darkmode.
\item
integration of multiple charts into one navigable collection.
\item
absolute node positioning for precise control.
\end{itemize}
\end{remark*}

\printbibliography

\end{document}